%Paper: alg-geom/9512007
%From: "Gang.XIAO" <xiao@aurora.unice.fr>
%Date: Sat, 9 Dec 1995 12:59:58 +0100

\magnification=1200
\voffset=0.3 true in
\parskip=2pt plus 1pt
\baselineskip=13.5pt plus 0.2pt minus 0.2pt
\lineskip=1.5pt plus0.5pt minus0.5pt
\lineskiplimit=0.5pt
\mathsurround=1pt
\font\tenmsx=msam10 \font\sevenmsx=msam7 \font\fivemsx=msam5
\newfam\msxfam
  \textfont\msxfam=\tenmsx \scriptfont\msxfam=\sevenmsx
\scriptscriptfont\msxfam=\fivemsx
\font\tenmsy=msbm10 \font\sevenmsy=msbm7 \font\fivemsy=msbm5
\newfam\msyfam
  \textfont\msyfam=\tenmsy \scriptfont\msyfam=\sevenmsy
\scriptscriptfont\msyfam=\fivemsy
\font\teneuf=eufm10 \font\seveneuf=eufm7 \font\fiveeuf=eufm5
\newfam\euffam
  \textfont\euffam=\teneuf \scriptfont\euffam=\seveneuf
\scriptscriptfont\euffam=\fiveeuf

\def\hexnumberA#1{\ifnum#1<10 \number#1\else
 \ifnum#1=10 A\else\ifnum#1=11 B\else\ifnum#1=12 C\else
 \ifnum#1=13 D\else\ifnum#1=14 E\else\ifnum#1=15 F\fi\fi\fi\fi\fi\fi\fi}
\def\MSX{\hexnumberA\msxfam}
\def\MSY{\hexnumberA\msyfam}
\def\EUF{\hexnumberA\euffam}

\def\text#1{\ifmmode
\mathchoice{\hbox{\tenrm #1}}
 {\hbox{\tenrm #1}}
 {\hbox{\sevenrm #1}}
 {\hbox{\fiverm #1}} \else\hbox{#1}\fi}
\def\matrix#1{\,\vcenter{\normalbaselines\mathsurround=0pt
  \ialign{\mathstrut\hfil$##$\hfil&&\quad\hfil$##$\hfil\crcr
  \mathstrut\crcr\noalign{\kern-\baselineskip}
  #1\crcr\mathstrut\crcr\noalign{\kern-\baselineskip}}}\,}
\def\fraction#1#2{{#1\over #2}}

\mathchardef\endproofmark="0\MSX03
\def\endofproof{\hphantom{\endproofmark}\hfill\llap{$\endproofmark$} }

\font\heading=cmr17 \font\footing=cmr9 \font\sectionfont=cmr10
\long\outer\def\Thm#1#2{\medbreak{\bf#1}{\sl#2}\medbreak}
\def\Prf#1{\medbreak{\bf#1}}

\def\Endrmk{\medbreak}
\def\Sct#1\par{\bigbreak\bigskip\bigskip\centerline{\sectionfont\S#1}
  \par\penalty10000\medskip}

\long\def\skipover#1\endskipover{}

\vglue1.5true in
\def\endproofmark{\hbox{\bf QED }}
\font\sectionfont=cmr12
\footline={\tenrm\hss --- \folio\ ---\hss}

\def\refindent{40pt}
\def\refhead{\sectionfont References}\def\reffont{\footing}
\def\Ref{\bigbreak\bigskip\bigskip\centerline{\refhead}\medskip
  \baselineskip=12pt plus.3pt minus.3pt\frenchspacing\reffont
  \everypar={\parindent=0pt\Refitem}\parindent=0pt
}

\def\Refitem [#1]#2:#3\par{\hangindent=\refindent\hangafter=1
  \rlap{[#1]}\hskip\hangindent #2: #3\par\smallbreak}
\medskip
{\heading
\par\centerline {            Non-symplectic involutions of a K3 surface}
}
\bigskip
{\footing
\par\centerline {                            Xiao, Gang}
\par\centerline {              D\'ept. de Math., UNSA, 06108 Nice, France}
\par\centerline {                    xiao@aurora.unice.fr}
}
\bigskip\bigskip
    Let $S$ be a smooth minimal $K3$ surface defined over $\mathchar"0\MSY43
$, $G$ a finite group
acting on $S$. The induced linear action of $G$ on $H^
{0}(\omega _{S})\cong \mathchar"0\MSY43 $ leads to an exact
sequence
                 $$1 \, \longrightarrow  \, K \, \longrightarrow
\, G \, \longrightarrow  \, N \, \longrightarrow
\, 1  \ ,$$
where the {\it non-symplectic part} $N$ is a cyclic group $\mathchar"0\MSY5A
_{m}$, which acts on the
intermediate quotient $S/K$ which is also $K3$. It is well-known that the Euler
number $\varphi (m)$ of $m$ must divide $22-\rho (S)$ ([N], Corollary 3.3), in
particular
$\varphi (m)\le 21$, hence $m\le 66$. It is also known that if $H$ is
non-trivial, then $S$ is
algebraic. In this case the quotient of $S$ by the action of $G$ is either an
Enriques surface or a rational surface. An example of $m=66$ has been
constructed in [K], where Kondo also gets the uniqueness of the $K3$ surface
with a non-symplectic action of $N\cong \mathchar"0\MSY5A
_{66}$, under the extra condition that $N$ acts
trivially on the N\'eron-Severi group of the surface. (Note that the
computation
in [K] contains an error, so that the case $m=44$ is missing in his final
result; the existence of this case is shown in our computation which follows.)

    The purpose of present article is to determine the $K3$ surfaces admitting
a
non-symplectic group $N$ of high order. More precisely, we look at the cases
                $$m=38, \, 44, \, 48, \, 50, \, 54,
\, 60, \text{ or }66.$$
\medskip
    \Thm{Theorem. }{ 1. There exists no $K3$ surface admitting a non-symplectic
$N$ of
order $60$.
\par     2. For each of the other 6 cases of $m$ as above, there is exactly one
$K3$
surface $S$ with $N\cong \mathchar"0\MSY5A _{m}$. The action of $N$ is also
unique (up to isomorphisms of $S$)
except in the case of $m=38$, in which case there are 2 different actions.}
\bigskip
    \Sct 1. General considerations\par
\medskip
    We consider the following situation: let $S$ be a $K3$ surface with a
non-symplectic automorphism group $G\cong \mathchar"0\MSY5A
_{m}$, i.e., no intermediate quotient of $S$ by
a subgroup of $G$ is again $K3$.
\par     Let $H\cong \mathchar"0\MSY5A _{t}$ be a subgroup of $G$, $X$ the
minimal resolution of singularities of
the intermediate quotient $S/H$, and let $\alpha \colon
\, \tilde S \, \longrightarrow  \, S$ be the minimal blow-up
such that the induced map $\pi \colon  \, \tilde S
\, \longrightarrow  \, X$ is a morphism. Let $B$ be the branch
locus of $\pi $. There is a $\mathchar"0\MSY51 $-divisor $\bar
\mathchar"0\EUF42 $ on $X$, supported on $B$, such that
$\alpha ^{*}(K_{S})\equiv \pi ^{*}(K_{X}+\bar \mathchar"0\EUF42
)$. If $B=\sum _{i}\Gamma _{i}$ is the decomposition of $B$ into irreducible
components, we have $\bar \mathchar"0\EUF42 ={\fraction
{1}{t}}\sum _{i}a_{i}\Gamma _{i}$, where the coefficient $a_
{i}$ is an integer with
$0\le a_{i}<t$ (cf. [X]).
\medskip
    \Thm{Lemma 1. }{ $B$ does not contain negative definite configurations of
$(-2)$-curves, therefore every component of $B$ has positive coefficient in
$\bar
\mathchar"0\EUF42 $.}
\par     \Prf{ Proof. }As $\pi ^{*}(K_{X}+\bar \mathchar"0\EUF42
)$ is nef, $K_{X}+\bar \mathchar"0\EUF42 $ is also nef. Therefore the
coefficients $a_{i}/t$ of components in a negative definite
$(-2)$-configuration
$\Gamma =\sum ^{k}_{i=1}\Gamma _{i}$ are equal to $0$. Then according to [X],
\S 1, $\Gamma
$ is the inverse image of
a singular point on $S/H$, as the coefficients $0$ are not of the form $1-1/n$
($n\ge 2$). This means that $\Gamma $ corresponds to an isolated fixed point
$p$ on $S$, for
the action of $H$. Furthermore if $K$ is the stabiliser of $p$, the
linearisation of
the action of $K$ on $T_{S}(p)$ is of the form $\left(\matrix
{\zeta  &0\cr 0 &\zeta ^{-1}\cr }\right)$, where $\zeta
$ is a root of unity
(cf. [BPV], \S III.5). This action being locally symplectic, the action of $K$
has
to be symplectic on $S$, which contradicts the hypothesis.
\par     For the second statement, we remark that [X], Lemma 4 is still true in
our case, so we can use [X], Lemma 5. \endofproof\Endrmk
\medskip
    \Thm{Lemma 2. }{ Let $G\cong \mathchar"0\MSY5A
_{m}$ be a group acting non-symplectically on a $K3$ surface
$S$. If $m>2$, the intermediate quotients of the action are all rational
surfaces.}
\par     \Prf{ Proof. }An intermediate quotient $X$ is an algebraic surface
with $p_
{g}=0$,
hence is either rational or Enriques. And a cyclic cover of $S$ over an
Enriques
surface must be non-ramified due to the above lemma, hence of degree 2 as the
$\pi _{1}$ of an Enriques surface is $\mathchar"0\MSY5A
_{2}$. Therefore:
\par     1. If $m$ is odd, all the intermediate quotients are rational.
\par     2. The quotient of a non-free action is rational.
\par     3. If $m=2n$ with $n$ odd, let $X$ be the intermediate quotient by
$\mathchar"0\MSY5A
_{n}$. Then the
quotient group $\mathchar"0\MSY5A _{2}$ acts on $X$, having a fixed point $p$.
The inverse image of $p$ on
$S$ has to contain a fixed point of the action of the subgroup
$\mathchar"0\MSY5A
_{2}$, as the order
of this inverse image is odd. Therefore the intermediate quotient of $S$ by
$\mathchar"0\MSY5A
_{2}$
is rational.
\par     4. If $m=4$, let $X$ be the intermediate quotient, $Y$ the final
quotient. If $X$
is Enriques, the quotient $\mathchar"0\MSY5A _{2}$-action on $X$ cannot have
fixed point, for
otherwise the inverse image of such a fixed point on $S$ has a
$\mathchar"0\MSY5A
_{2}$-stabiliser
different than the first $\mathchar"0\MSY5A _{2}$-subgroup, which implies
$G\cong
\mathchar"0\MSY5A _{2}^{2}$, impossible. However
an Enriques surface does not allow fixed-point free involutions, as e.g.
$\chi (\mathchar"024F _{X})=1$ is not divisible by 2.
\par     Now in the general case, a $\mathchar"0\MSY5A
_{2}$-subgroup of $\mathchar"0\MSY5A _{m}$ is contained in a subgroup
$\mathchar"0\MSY5A
_{k}$
with either $k=4$ or $k=2n$ where $n$ is odd. This proves the lemma because any
quotient of a rational surface is rational. \endofproof\Endrmk
\medskip
    Similarly, one shows

    \Thm{Lemma 3. }{ Let $G\cong \mathchar"0\MSY5A
_{n^{2}}$ acting non-symplectically on a $K3$ surface $S$ where $n$
is a prime, and let $H\cong \mathchar"0\MSY5A _{n}$ be the subgroup of $G$,
$Q=G/H$, $X=S/H$. Let $D$ be the
branch locus of the projection $S \, \longrightarrow
\, X$. Then all the fixed points of the
induced action of $Q$ on $X$ are located on $D$.}
\par     \Prf{ Proof. }Let $p$ be such a fixed point. If it is not on $D$, its
inverse image
on $S$ is composed of $n$ points, therefore each of them has a stabiliser
isomorphic to $\mathchar"0\MSY5A _{n}$ in $G$, different from $H$. This is
impossible as $G$ is cyclic. \endofproof\Endrmk
\medskip
    Now let $S$ be a $K3$ surface with a non-symplectic action of
$G=\mathchar"0\MSY5A
_{m}$ where $m>2$
is even. Let $X$ be the intermediate quotient of $S$ by the unique
$\mathchar"0\MSY5A
_{2}$-subgroup $<\iota >$
of $G$. $X$ is a smooth rational surface. Let $B$ be the branch locus of the
projection $\pi \colon  \, S \, \longrightarrow  \,
X$. $B$ is a smooth divisor linearly equivalent to $-2K_
{X}$. We
have
                 $$10-K_{X}^{2} \, = \, \rho (X) \,
\le  \, \rho (S) \, \le  \, 22-\varphi (m)  \ .$$
Let $Q$ be the quotient of $G$ by $\mathchar"0\MSY5A
_{2}$, which acts naturally on $X$. $B$ is invariant
under this action.
\medskip
    \Thm{Lemma 4. }{ If $X\cong \mathchar"0\MSY50
^{2}$, then either $m\le 30$, $m=42$, or $m=50$.}
\par     \Prf{ Proof. }$B$ is a smooth sextic.
\par     Note first that an action of $\mathchar"0\MSY5A
_{2}$ on $X$ always has a fixed point plus a fixed
line, hence by Lemma 3, $m/2$ must be odd.
\par     Let $\gamma $ be a generator of $Q$. The action of $\gamma
$ on $X$ has either a fixed
point $p$ and a line $L$ composed of fixed points; or 3 fixed points $p_
{1},p_{2},p_{3}$.
\par     In the first case, let $H$ be a general line passing through $p$. $H$
is
invariant, and the action of $Q$ on $H$ has exactly 2 fixed points, namely $p$
and
$H\cap L$. But then the intersection $H\cap B$ has to be invariant; as $|H\cap
B|=6$ and $Q$ is
cyclic, we must have $|Q|\le 5$.
\par     For the second case, assume first that $B$ meets each line $L_
{i}$ passing
through $p_{i}$ and $p_{i+1}$ (letting $p_{4}=p_{1
}$) only on $p_{i}$ and $p_{i+1}$. By the smoothness of
$B$, this is possible only when, say, $B$ is tangent to $L_
{i}$ to order $5$ at $p_{i}$ for
$i=1,2,3$. Consider the projection $f\colon  \, B
\, \longrightarrow  \, B/Q=C$. It is clear that $f$ is
ramified exactly at the 3 points $p_{i}$, hence by Hurwitz Formula, one gets
$|Q|=3,7$ or $21$.
\par     Finally, assume that $B\cap L_{1}$ contains a point other than $p_
{1}$ and $p_{2}$. Because
the set $B\cap L_{1}$ is invariant under the action of $Q$, The subgroup $H$ of
$Q$ fixing
every point of $L_{1}$ is of index at most 5. Also $|H|\le
5$ as in the first case, and
we get the conclusion of the lemma. \endofproof\Endrmk
\medskip
    Now assuming $\rho (X)>1$, we have ``ruling''s on $X$, i.e., a morphism
$r\colon
\, X \, \longrightarrow  \, $
$C\cong \mathchar"0\MSY50 ^{1}$ whose general fibres are isomorphic to
$\mathchar"0\MSY50
^{1}$. The pull-back of $r$ on $S$ is an
elliptic fibration.
\par     By Hurwitz Formula, the induced cover $r|_
{B}\colon  \, B \, \longrightarrow  \, C$ has total ramification
index $\delta \le 24$.
\medskip
    \Thm{Lemma 5. }{ Let $\sigma $ be a non-symplectic automorphism in $Q$
which fixes each
fibre of a ruling $r\colon  \, X \, \longrightarrow
\, C$. $\sigma $ is either trivial or isomorphic to $\mathchar"0\MSY5A
_{3}$. In the
latter case $B$ contains a section $C_{0}$ of $r$ with $C_
{0}^{2}=-4$.}
\par     \Prf{ Proof. }Let $K$ be the inverse image of $<\sigma
>$ in $G$. $K$ acts on the inverse
image $E$ of a general fibre $F$ of $r$, which is an elliptic curve. As $K$ is
cyclic
and contains the elliptic involution, one must have $K=\mathchar"0\MSY5A
_{2}, \, \mathchar"0\MSY5A _{4}$ or $\mathchar"0\MSY5A
_{6}$.
\par     Moreover in the case of $\mathchar"0\MSY5A
_{4}$, the two fixed points of $\sigma $ on $F$ must be
contained in $B$. This implies a decomposition $B=B_
{1}+B_{2}$, with $B_{1}$ and $B_{2}$ both of
degree 2 over $C$, and $B_{1}B_{2}=0$. As $K_{X}^{2
}\ge 6$, one sees easily that this cannot
happen, say, by contracting $X$ into a Hirzebruch surface.
\par     In the case of $\mathchar"0\MSY5A _{6}$, the existence of $C_
{0}$ is due to the existence of a total
fixed point for the action of $K$ on $E$; and $C_{0
}^{2}=-4$ is dictated by the condition
$B\equiv -2K_{X}$. \endofproof\Endrmk
\bigskip
    \Prf{ Definition. }Let $Y=\mathchar"0\MSY46 _{e
}$ be a Hirzebruch surface of invariant $e$ with the
ruling $r\colon  \, Y \, \longrightarrow  \, C\cong
\mathchar"0\MSY50 ^{1}$, and let $\gamma $ be an automorphism of finite order
$n$ on $Y$
respecting $r$, such that its induced action on $C$ is also of order $n$. Let
$F_
{1},F_{2}$
be the two invariant fibres of $r$.
\par     For any fixed point $p$ of $\gamma $, define the {\it type} of $p$,
$\tau
_{p}$, as follows.
Choose local parameters $\left\lbrace t,x\right\rbrace
$ of $p$, where $x$ is vertical with respect to $r$,
such that the action of $\gamma $ diagonalizes: $\gamma
(t)=\xi t$, $\gamma (x)=\xi ^{\alpha }x$, where $\xi
$ is a
primitive $n$-th root of unity, $0\le \alpha <n$. And define $\tau
_{p}=\alpha $. Note that $\tau _{p}$ depends
only on the action of the group $<\gamma >$.
\par     When $e>0$, let $C_{0}$ be the section of negative self-intersection
on $Y$; when
$e=0$, we fix an invariant flat section to be $C_{0
}$. With respect to $C_{0}$, we may
define the {\it type} of $F_{i}$, $\tau _{i}$, to be $\tau
_{F_{i}\cap C_{0}}$.
\par     Note that if $p$ and $q$ are two fixed points on a same fibre $F_
{i}$, we have
                $$\tau _{p}+\tau _{q}\equiv 0 \, \pmod
{n} \, .$$\Endrmk
\medskip
    \Thm{Lemma 6. }{ $\tau _{1}+\tau _{2}+e\equiv
0 \, \pmod{n}$.}
\par     \Prf{ Proof. }Let $p_{i}=F_{i}\cap C_{0}$, and let $Y'$ be the surface
resulting from an
elementary transform centered at $p_{1}$. As $p_{1
}$ is fixed under $\gamma $, we have an
induced action on $Y'$, for which the type of $F_{1
}$ becomes $\tau _{1}-1$. This allows us to
show the lemma only for the case $\tau _{1}=\tau _
{2}=0$, but in this case $\gamma $ has no isolated
fixed point, hence the quotient $Y/<\gamma >$ is smooth Hirzebruch surface
$\mathchar"0\MSY46
_{d}$, so that
$e=nd$. \endofproof\Endrmk
\medskip
    \Thm{Lemma 7. }{ Let $X$ be a smooth rational surface with $K_
{X}^{2}>0$, and let
                 $$|F_{1}| \, ,\ldots , \, |F_{n}|$$
be $n$ rulings with $F_{i}F_{j}=a$, $\forall i,j$. Then
                  $$K_{X}^{2}\le {\fraction {  4n}
{a(n-1)}} \, .$$}
\par     \Prf{ Proof. }Let $D={\fraction {  2}{a(n-1)
}}\left(F_{1}+F_{2}+\cdots +F_{n-1}\right)$. As $(K_
{X}+D)F_{n}=0$, we have
                $$K_{X}^{2}-{\fraction {  4n}{a(n-1)
}}=\left(K_{X}+D\right)^{2}\le 0$$
by Hodge Index Theorem. \endofproof\Endrmk
\medskip
    \Thm{Lemma 8. }{ In the case where $\rho (X)>1$ and $m=38$ or $m\ge
44$, $X$ has an equivariant
ruling under the action of $Q$.
\par     Moreover, when $3|m$, the ruling is invariant under the subgroup of
order 3.}
\par     \Prf{ Proof. }When $3{\not|}m$ (hence $K_
{X}^{2}\ge 6$) or $\varphi (m)=20$, the above Lemma 7 tells that
the orbit of a ruling under $Q$ has at most 2 elements, with fibres
intersecting
each other by 1. Hence the only possibility to exclude is that $X$ contracts to
a $X_{0}\cong \mathchar"0\MSY50 ^{1}\times \mathchar"0\MSY50
^{1}$, with the action of $Q$ exchanging the two factors. As $|Q|$ is not
divisible by 4, the subgroup $H$ of order 2 of $Q$ acts on $X_
{0}$ by exchanging the
factors. But then all the points on the diagonal $D$ are fixed under $H$, hence
$D$
is contained in the image $B_{0}$ of $B$, but then $D(B_
{0}-B)=6$, and we cannot blow up
$X_{0}$ at most 2 times to make $B$ smooth.
\par     Therefore we can assume that there is an element $\sigma
$ of order 3 in $Q$. We
first show that there is an equivariant ruling under $<\sigma
>$. To do so let $|F_{1}|, \, $
$|F_{2}|, \, |F_{3}|$ be 3 rulings forming an orbit of $<\sigma
>$. Lemma 7 forces $F_{i}F_{j}=1$ for
$i\ne j$, hence there exists a contraction $v\colon
\, X \, \longrightarrow  \, X_{0}\cong \mathchar"0\MSY50
^{2}$ such that the images of
the pencil $|F_{i}|$ is a pencil of lines through a point $p_
{i}$, for $i=1,2,3$. The
contraction $v$ is unique when the points $p_{i}$ are colinear; and there is
exactly
one other such contraction when the points are not colinear. In any case,
there is a subgroup $H$ of index $\le 2$ in $G$ which has an induced action on
$X_
{0}$.
\par     Note that the action of $\sigma $ on $X_{0
}$ cannot fix a singular point of $B_{0}=v(B)$,
for otherwise the pull-back of the pencil of lines through such a point would
give rise to an equivariant ruling for $<\sigma >$. Therefore the number of
singular
points of $B_{0}$ is divisible by 3. As this number is at most 5, $B_
{0}$ has to be
smooth outside $\left\lbrace p_{1},p_{2},p_{3}\right\rbrace
$, and $K_{X}^{2}=6$.
\par     Let $K\subset H$ be the stabiliser of $p_
{1}$. As $K$ fixes also $p_{2},p_{3}$ as well as at least
3 fixed points of the action of $\sigma $ on $X_{0
}$, the only way for $K$ to have a
non-trivial action on $X_{0}$ is that $p_{1},p_{2},p_
{3}$ are on a same line $L$ which is then
fixed pointwise by $K$. As $B_{0}$ has either ordinary double point or ordinary
cusp
on $p_{i}$ and $|K|>2$, the local invariance of $B_
{0}$ around $p_{i}$ forces $L$ to be a
component of $B_{0}$, which is impossible as $B_{0
}(B_{0}-L)=5>3$.
\par     So now we have a ruling $r\colon  \, X \,
\longrightarrow  \, C$ which is equivariant under $\sigma
$. When $r$
is invariant, it is easy to see that it is equivariant under $Q$: indeed, let
$p$
be a general point in $X$, $\Sigma $ the orbit of $p$ under $<\sigma
>$, $F$ the fibre containing $p$,
and let $\gamma \in Q$. By the commutativity of $Q$, $\gamma
$ sends $\Sigma $ to an orbit $\Sigma '$ of $<\sigma
>$, which
is contained in a fibre $F'$ of $r$. Now if $\gamma
(F)\ne F'$, we would have $\gamma (F)F\ge 3$ as
$\Sigma '\subseteq F'\cap \gamma (F)$, which contradicts Lemma 7 (by taking
$n=2$).
\par     It remains to exclude the case where $r$ is equivariant but not
invariant
under $\sigma $. Let $\tilde r\colon  \, S \, \longrightarrow
\, C$ be the pull-back of $r$ on $S$, $\tilde \sigma
$ the element of order
3 in $G$ whose image in $Q$ is $\sigma $. In this case the fixed locus of
$\tilde
\sigma $ is contained
in two fibres of $\tilde r$, hence is composed of $e_
{1}$ isolated fixed points, $e_{2}$
rational curves of self-intersection $-2$, plus possibly one or two elliptic
curves. Let $\alpha \colon  \, \hat S \, \longrightarrow
\, S$ be the blow-up of the isolated fixed points of $\tilde
\sigma $.
Then the quotient $Y=\hat S/<\tilde \sigma >$ is a smooth rational surface with
$K_
{Y}^{2}=-(e_{1}+8e_{2})/3$,
$c_{2}(Y)=8+(5e_{1}+4e_{2})/3$. Hence $e_{1}-e_{2}=3$ as $K_
{Y}^{2}+c_{2}(Y)=12$, but then
        $$\rho (S)=\rho (\hat S)-e_{1}\ge \rho (Y)-e_
{1}=10+(-2e_{1}+8e_{2})/3=8+2e_{2}\ge 8$$
which is excluded by our conditions. \endofproof\Endrmk
\bigskip\bigskip
    The following remark is useful for the existence of the cases.
\medskip
    \Thm{Lemma 9. }{ An automorphism $\gamma $ on $X$ lifts up to an
automorphism on $S$ if and
only if $\gamma (B)=B$.}
\par     \Prf{ Proof. }The double cover $\pi $: $S
\, \longrightarrow  \, X$ is determined by an element $\delta
\in Pic(X)$
such that $B\equiv 2\delta $. As $X$ is simply connected, $\delta
$ hence $\pi $ is determined by $B$. \endofproof\Endrmk
\bigskip\bigskip
    \Sct 2. The cases with $3|m$\par
\bigskip
    We consider in this section the cases $m=48, \,
54, \, 60, \, 66$. According to
Lemma~8, we have a ruling $r\colon  \, X \, \longrightarrow
\, C$ which is equivariant under $Q$, and such
that the action of the subgroup $<\sigma >$ of order 3 on $X$ has a fixed locus
composed
of two sections $C_{0},C_{1}$, one of which, say $C_
{0}$, is a component of $B$.
\par     There is a unique contraction $t_{1}\colon
\, X \, \longrightarrow  \, X_{1}$ to a Hirzebruch surface $X_
{1}$
with respect to $r$, such that the image of $C_{0}$ is still of
self-intersection $-4$.
The action of $\sigma $ descends to $X_{1}$, with projection $t_
{2}\colon  \, X_{1} \, \longrightarrow  \, X_{2}=X_
{1}/<\sigma >$, where
$X_{2}\cong \mathchar"0\MSY46 _{12}$, and a ruling $r_
{2}\colon  \, X_{2} \, \longrightarrow  \, C$ induced from $r$.
\par     We have 3 sections $C_{2},C_{3},C_{4}$ of $r_
{2}$, with $-C_{2}^{2}=C_{3}^{3}=C_{4}^{2}=12$, such that $C_
{2}+C_{3}$
is the branch locus of $t_{2}$, and $C_{2}+C_{4}$ is the image of $B$. There is
an induced
action of $\bar Q=Q/<\sigma >\cong \mathchar"0\MSY5A
_{m/6}$ on $X_{2}$, respecting $r_{2}$. Let $F_{1},F_
{2}$ be the two invariant
fibres of $r_{2}$ under this action, and let $\alpha
_{i}$ be the number of intersection of
$C_{3}$ and $C_{4}$ on $F_{i}$. Because $C_{3}C_{4
}=12$, we have clearly $\alpha _{1}+\alpha _{2}=12-m/6$. Assume $\alpha
_{1}\le \alpha _{2}$.
\par     Let $\tau _{i}$ be as in the definition above Lemma 6, for the action
of $\bar
Q$ on $X_{2}$.
We have $\tau _{i}=m/6-\alpha _{i}$ as $C_{2},C_{3
},C_{4}$ are invariant curves. Let $p_{i}=C_{3}\cap
F_{i}$, $q_{i}=C_{2}\cap F_{i}$.
As in the proof of Lemma 6, after $\alpha _{1}$ successive elementary
transformations
centered on $p_{1}$ and $\alpha _{2}$ transformations centered on $p_
{2}$, we get a surface $X_{3}\cong \mathchar"0\MSY46
_{m/6}$
on which $\bar Q$ acts without isolated fixed point; Therefore the quotient $X_
{4}=X_{3}/\bar Q$
is the Hirzebruch surface $\mathchar"0\MSY46 _{1}$. Contracting the negative
section of $X_
{4}$, we
arrive at the projective plane on which the images of the ramification curves
$C_{3},C_{4},F_{1},F_{2}$ form four lines with normal crossings. Such a
configuration being
unique, the uniqueness of $S$ for each $m$ will be shown if we can show the
uniqueness of the couple $(\alpha _{1},\alpha _{2})$ for each $m$.
\medskip
    For $m=66$, the unique possibility is $(\alpha
_{1},\alpha _{2})=(0,1)$; for $m=60$, $(\alpha _{1
},\alpha _{2})=(0,2)$
or $(1,1)$. $(0,2)$ is impossible because the subgroup of order 2 in $Q$ would
contradict Lemma 3, as (the strict transform on $X$ of) $F_
{1}$ is clearly not in $B$.
While in the case of $(1,1)$, let $\tilde \gamma $ be an element of order 5 in
$G$, $\gamma
$ the image
of $\tilde \gamma $ in $Q$. The action of $\gamma
$ on $T_{X_{2}}(q_{1})$ is by definition of the form $\left(\matrix
{\zeta  &0\cr 0 &\zeta ^{-1}\cr }\right)$
where $\zeta $ is a root of unity of order 5; but then the action of $\tilde
\gamma $ on the
inverse image of $q_{1}$ is also of the form $\left(\matrix
{\zeta  &0\cr 0 &\zeta ^{-1}\cr }\right)$ because $6\equiv
1\pmod5$, which
means that $\tilde \gamma $ is a symplectic automorphism. This shows the
non-existence of
$m=60$.
\par     For the same reason, the case $m=54$ admits only $(\alpha
_{1},\alpha _{2})=(1,2)$ because $(0,3)$
does not verify Lemma 3 with respect to the subgroup of order 3 in $Q$. And the
case $m=48$ admits only $(\alpha _{1},\alpha _{2})=(1,3)$ by considering the
subgroup of order 2 in
$Q$.
\medskip
    Finally, the existence of the cases $48,54,66$ can be shown by reversing
the
above argument: take 2 fibres $F'_{1},F'_{2}$ as well as 3 sections $C'_
{2},C'_{3},C'_{4}$ on the
Hirzebruch surface $\mathchar"0\MSY46 _{1}$, with $-C_
{2}'^{2}=C_{3}'^{2}=C_{4}'^{2}=1$. Make a cyclic cover $X_
{3} \, \longrightarrow  \, \mathchar"0\MSY46 _{1}$
of order $m/6$ ramified along $F'_{1}$ and $F'_{2}$, and note by $F_
{1}$, etc. the inverse
image of $F'_{1}$, etc. Make $\alpha _{i}$ elementary transforms on $q_
{i}=F_{i}\cap C_{2}$ for $i=1,2$ to get
the surface $X_{2}$, then a triple cover $t_{2}\colon
\, X_{1} \, \longrightarrow  \, X_{2}$ ramified along $C_
{2}$ and $C_{3}$,
and blow up the singularities of the inverse image of $C_
{4}$ to get $t_{1}\colon  \, X \, \longrightarrow
\, X_{1}$.
It is easy to see that the map $X \, \mathrel{\smash-\smash-\mathchar"0221
} \, \mathchar"0\MSY46 _{1}$ thus constructed is generically
cyclic of order $m/2$, and we can use Lemma 9 to see that this cyclic action of
order $m/2$ on $X$ lifts to an automorphism group $G$ of order $m$ on $S$. It
remains
only to verify that $G$ acts non-symplectically, for which it suffices to
verify
that every minimal subgroup of $G$ acts non-symplectically, which can be done
locally around a fixed point. Details of the verification are left to the
reader.
\bigskip
    \Sct 3. The remaining cases\par
\bigskip
    The case $m=50$:
\par     We have shown in \S 1 that $X\cong \mathchar"0\MSY50
^{2}$, and that the action of $Q=<\gamma >$ is of the form
$\gamma (x_{0}:x_{1}:x_{2})=(\zeta x_{0}:\zeta ^{5\alpha
+1}x_{1}:x_{2})$, where $\zeta $ is a primitive root of unity of order
$25$, and $\alpha \in \mathchar"0\MSY5A $. Letting $p_
{1}=(1:0:0)$, $p_{2}=(0:1:0)$, $p_{3}=(0:0:1)$, $B$ intersects
$L_{1}-\left\lbrace p_{1},p_{2}\right\rbrace $ transversally at 5 points, hence
it passes through, say, $p_
{2}$. As $B$
cannot intersect $L_{2}$ and $L_{3}$ at points other than $p_
{1},p_{2},p_{3}$, we must have
$B\cap L_{3}=\left\lbrace p_{3}\right\rbrace $ with a tangent of order 6.
Therefore a local computation at $p_
{3}$ gives
$\alpha =1$. Also the intersection of $B$ with $L_
{3}$ shows that the equation of $B$ contains
the term $X_{0}^{6}$, with $\gamma (X_{0}^{6})=\zeta
^{6}X_{0}^{6}$. There are only two other monomials of degree
6 with the same character, namely $X_{0}X_{1}^{5}$ and $X_
{1}X_{2}^{5}$. One concludes easily that
modulo automorphisms of $X$, the equation of $B$ is
                    $$X_{0}^{6}+X_{0}X_{1}^{5}+X_{1
}X_{2}^{5}=0  \ .$$
This proves the uniqueness as well as the existence in view of Lemma 9.
\par     Passing to the total quotient, one sees that $S$ is the smooth minimal
model
of a cyclic cover of $\mathchar"0\MSY50 ^{2}$ ramified along 4 lines of general
position, with
respective ramification indices 2, 5, 25, 50.
\medskip
    The case $m=44$:
\par     Let $F_{1},F_{2}$ be the two invariant fibres of $r\colon
\, X \, \longrightarrow  \, C$ under the action of
$Q$. $r|_{B}$ has two ramifications on $F_{1}+F_{2
}$.
\par     Note that if $r|_{B}$ has at most one ramification on a fibre $F_
{i}$, then $B\cap F_{i}$ has
at least 3 points, so $\tau _{i}=0$ for the action of the subgroup
$\mathchar"0\MSY5A
_{11}$ of $Q$. This
excludes the case where the two ramifications are distributed on the two
invariant fibres, as in this case $\tau _{1}=\tau
_{2}=e=0$ for $\mathchar"0\MSY5A _{11}$, which is impossible
because the horizontal degree of $B$ is not a multiple of 11.
\par     We may thus assume that $B$ is tangent to $F_
{1}$ of order 3 at a point $p_{1}$. Then
$11|\tau _{2}$ and $11{\not|}\tau _{1}$ for the action of $Q$, so $e>0$. In
fact the local invariance of $B$
at $p_{1}$ gives $\tau _{p_{1}}=15$, and Lemma 6 gives quickly $\tau
_{1}=7$, $\tau _{2}=11$, $e=4$, and then a
disjoint decomposition $B=B_{0}+C_{0}$ with $B_{0}$ smooth irreducible.
\par     After 7 successive elementary transforms centered at $F_
{1}\cap C_{0}$ then 11
elementary transforms centered at the fixed point of $F_
{2}$ not on $C_{0}$, we get a
surface $X_{1}\cong \mathchar"0\MSY46 _{0}$. Let $X_
{2}\cong \mathchar"0\MSY46 _{0}$ be its quotient by $Q$, and let $B_
{2},C_{2},F_{3},F_{4}$ be
respectively the images on $X_{2}$ of $B_{0},C_{0},F_
{1},F_{2}$. $B_{2}$ is smooth of bidegree $(3,1)$,
totally tangent to $F_{3}$ and tangent to $F_{4}$ of order 2 at the point where
the
horizontal section $C_{2}$ passes through. Such a configuration being unique up
to
automorphisms of $\mathchar"0\MSY46 _{0}$, we get the uniqueness of this case.
And the existence is
shown by reversing the arguments, as for the preceding cases. (To see that the
action is non-symplectic, just note that as there is no symplectic
automorphism of order 11, one has only to show that there is a cyclic subgroup
of order 4; this can be done locally around a fixed point.)
\par     $S$ is birationally a cyclic cover of $\mathchar"0\MSY50
^{1}\times \mathchar"0\MSY50 ^{1}$ ramified along $B_
{2},C_{2},F_{3},F_{4}$, with
respective ramification indices $2, \, 2, \, 44, \,
44$.
\medskip
    The case $m=38$:
\par     Choose a contraction $\sigma \colon  \, X
\, \longrightarrow  \, X_{0}\cong \mathchar"0\MSY46
_{e}$ onto a Hirzebruch surface
$r_{0}\colon  \, X_{0} \, \longrightarrow  \, C$, and let $B_
{0}$ be the image of $B$ on $X_{0}$, and $F_{1},F_
{2}$ the invariant
fibres of $r_{0}$. Let $\beta _{i}$ be the number of ramifications of $r_
{0}|_{B_{0}}$ on $F_{i}$. We have
$\beta _{1}+\beta _{2}=5$, and can assume $\beta _
{1}<\beta _{2}$.
\par     For any fixed point $p$ of the action of $Q$ on $X_
{0}$, we have $\tau _{p}>1$: indeed,
otherwise as $e\le 4$, after at most 6 elementary transforms, we get a surface
$X'\cong \mathchar"0\MSY50 ^{1}\times \mathchar"0\MSY50
^{1}$, such that the induced action of $Q$ fixes one fibre pointwise. But
then Lemma 6 says that it is the pull-back of an action on $\mathchar"0\MSY50
^{1}$, hence the
strict transform $B'$ of $B_{0}$ on $X'$ should have a horizontal degree
divisible by
$19$, or $B'^{2}\ge 152$. This is impossible because $B_
{0}^{2}=32$ and each elementary
transform increases the square by at most $16$.
\par     One sees from this remark that $B_{0}$ meets each $F_
{i}$ at at most 2 points,
and that if $B_{0}$ have an ordinary double point, then one of the branches is
tangent to the fibre. And a local computation of $\tau
$ shows that $B_{0}$ cannot be
tangent to $F_{1}$ at two points. Therefore $\beta
_{1}=2$, and there is a point $p_{1}$ at which
$B_{0}$ is tangent to $F_{1}$ of order 3, with $\tau
_{p_{1}}=13$. $B_{0}\cap F_{1}$ contains another point $q_
{1}$
of transversal intersection.
\par     Now that $\beta _{2}=3$, one sees quickly that there are only two
possibilities
satisfying the above conditions: either $B_{0}\cap
F_{2}$ contains one point $p_{2}$ which is
tangent of order 4, or $B_{0}\cap F_{2}=\left\lbrace
p_{2},q_{2}\right\rbrace $ where $p_{2}$ is an ordinary double point of
$B_{0}$ with one branch tangent to $F_{2}$.
\medskip
    In the first possibility, $\tau _{p_{2}}=5$ and Lemma 6 leaves only one
possibility
$\tau _{1}=13$, $\tau _{2}=5$, $e=1$, with the negative section $C_
{0}$ passing through $p_{1}$ and $p_{2}$.
\par     After 6 successive elementary transforms centered on $q_
{1}$ and 5 on $p_{2}$ then
passing to quotient of $Q$, we get a $X_{1}\cong \mathchar"0\MSY50
^{1}\times \mathchar"0\MSY50 ^{1}$, with the image $B_
{1}$ of $B_{0}$ which is
smooth of bidegree $(4,1)$, intersecting $F_{3}$ at two points with one
transversal;
and tangent to $F_{4}$ at one point of order 4, where $F_
{3},F_{4}$ are respectively the
images of $F_{1},F_{2}$. Such a configuration being unique (it is the graph of
a map
$f\colon  \, \mathchar"0\MSY50 ^{1} \, \longrightarrow
\, \mathchar"0\MSY50 ^{1}$ determined by a pencil generated by two divisors
$4s_
{1}$ and $3s_{2}+s_{3}$,
hence is unique modulo automorphisms of the first $\mathchar"0\MSY50
^{1}$), we get the uniqueness
as well as the existence of this case:
\par     $S$ is birationally a cyclic cover of $\mathchar"0\MSY50
^{1}\times \mathchar"0\MSY50 ^{1}$ ramified over $B_
{1},F_{3}$ and $F_{4}$, with
respective ramification indices $2, \, 19, \, 38$.
\medskip
    In the second possibility, $\tau _{p_{2}}=10$ so $\tau
_{q_{2}}=9$. And we can choose the
contraction $\sigma $ such that $e=4$, and $q_{1},q_
{2}$ are on the negative section $C_{0}$. This
gives rise to a disjoint decomposition $B_{0}=B'_{0
}+C_{0}$, and after elementary
transforms centered on $p_{1}$ and $q_{2}$ then passing to the quotient, we get
a
$X_{1}\cong \mathchar"0\MSY50 ^{1}\times \mathchar"0\MSY50
^{1}$ with a same configuration as in the case $m=44$, hence the uniqueness
and the existence of this case.
\medskip
    \Prf{ Remark. }It is easy to see that the $K3$ surface $S$ in the two cases
of $m=38$
are the same, by analysing the elliptic fibration induced by $r$. The two
different actions arise from the choice of the involution.\Endrmk
\medskip
\Ref

    [BPV] Barth, W. / Peters, C. / Van de Ven, A.: Compact Complex Surfaces,
Springer 1984

    [K] Kondo, S.: On automorphisms of algebraic $K3$ surfaces which act
trivially on Picard groups, Proc. Japan Acad. 62 (A), 356-359 (1986)

    [N] Nikulin, V. V.: Finite automorphism groups of K\"ahler $K3$
surfaces, Trans. Moscow Math. Soc. 38, 71-137 (1980)

    [X] Xiao, G.: Bound of automorphisms of surfaces of general type I,
Annals of Math. 139, 51-77 (1994)
\medskip
\par\end